\documentclass[11pt]{article}

\usepackage[dvips]{graphicx}
\usepackage{amsfonts}
\usepackage{amssymb}
\usepackage{amsthm}
\usepackage{epsfig}

\theoremstyle{plain}
\newtheorem{theorem}{Theorem}[section]
\newtheorem{lemma}[theorem]{Lemma}
\newtheorem{proposition}[theorem]{Proposition}
\newtheorem{corollary}[theorem]{Corollary}
\theoremstyle{definition}
\newtheorem{definition}[theorem]{Definition}
\newtheorem{discussion}[theorem]{Discussion}
\newtheorem{remark}[theorem]{Remark}

\newtheorem{example}[theorem]{Example}

\newcommand{\todo}[1]{\vspace{5 mm}\par \noindent      
\marginpar{\textsc{ToDo}}
\framebox{\begin{minipage}[c]{0.95 \textwidth}
\tt #1 \end{minipage}}\vspace{5 mm}\par}

\newcommand{\idiot}[1]{\vspace{5 mm}\par \noindent     
\framebox{\begin{minipage}[c]{0.95 \textwidth}
\tt #1 \end{minipage}}\vspace{5 mm}\par}

\renewcommand{\todo}[1]{}
\renewcommand{\idiot}[1]{}

\newcommand{\xs}{x_1,\ldots,x_n}                
\newcommand{\vs}{v_1,\ldots,v_n}                
\newcommand{\Ms}{M_1,\ldots,M_q}                
\newcommand{\Fs}{F_1,\ldots,F_q}                

\newcommand{\m}{{\mathbf{m}}\ }                 
\newcommand{\height}{{\rm{height}}\ }           
\newcommand{\dimn}{ {\rm{dim}} \ }
\newcommand{\F}{{\mathcal{F}}}                  
\newcommand{\N}{{\mathcal{N}}}                  
\newcommand{\U}{{\mathcal{U}}}                  
\newcommand{\R}{{\mathcal{R}}}                  
\newcommand{\D}{\Delta}                         
\newcommand{\al}{\alpha}                         
\newcommand{\be}{\beta}                          
\newcommand{\df}{\delta_{\mathcal{F}}}          
\newcommand{\dn}{\delta_{\mathcal{N}}}          
\newcommand{\ndiv}{\not |}              
\newcommand{\st}{\ | \ }                        
\newcommand{\tuple}[1]{\langle #1 \rangle}      
\newcommand{\rmv}[1]{\setminus \langle #1\rangle} 
\newcommand{\void}[1]{}

\title{\sc Cohen-Macaulay Properties of Square-Free Monomial Ideals}
\author{Sara Faridi\thanks{University of Ottawa, 585 King Edward Ave,
Ottawa, ON K1N 6N5, Canada.  email: \emph{faridi@ uottawa.ca.}\newline
2000 Mathematics Subject classification: 13, 05. \newline This
research was supported by an NSERC Postdoctoral Fellowship.}}

\begin{document}

\maketitle

\begin{abstract} In this paper we study simplicial complexes as higher 
  dimensional graphs in order to produce algebraic statements about
  their facet ideals.  We introduce a large class of square-free
  monomial ideals with Cohen-Macaulay quotients, and a criterion for
  the Cohen-Macaulayness of facet ideals of simplicial trees. Along
  the way, we generalize several concepts from graph theory to
  simplicial complexes.
\end{abstract}

\section{Introduction}
From the point of view of commutative algebra, the focus of this paper
is on finding square-free monomial ideals that have Cohen-Macaulay
quotients. In~\cite{Vi1} Villarreal proved a criterion for the
Cohen-Macaulayness of edge ideals of graphs that are trees. Edge
ideals are square-free monomial ideals where each generator is a
product of two distinct variables of a polynomial ring. These ideals
have been studied extensively by Villarreal, Vasconcelos and Simis
among others. In~\cite{Fa} we studied a generalization of this
concept; namely the \emph{facet ideal} of a simplicial complex. By
generalizing the definition of a ``tree'' to simplicial complexes, we
extended the results of~\cite{SVV} from the class of edge ideals to
all square-free monomial ideals.

Below we investigate the structure of simplicial complexes in order to
show that Villarreal's Cohen-Macaulay criterion for graph-trees
extends to simplicial trees (Corollary~\ref{CM-criterion}). This is of
algebraic and computational significance, as it provides an effective
criterion for Cohen-Macaulayness that works for a large class of
square-free monomial ideals. We introduce a condition on a simplicial
complex that ensures the Cohen-Macaulayness of its facet ideal, and a
method to build a Cohen-Macaulay ideal from any given square-free
monomial ideal. Along the road to the algebraic goal, this study sheds
light on the beautiful combinatorial structure of simplicial
complexes.

The paper is organized as follows: Sections~\ref{basic-setup}
to~\ref{basic-trees-section} review the basic definitions and cover
the elementary properties of trees.  In
Section~\ref{n-partite-section} we draw comparisons between graph
theory and simplicial complex theory, and prove a generalized version
of K\"{o}nig's Theorem in graph theory for simplicial complexes.  We
then go on to prove a structure theorem for unmixed trees in
Section~\ref{structure-theorem-section}. We introduce the notion of a
\emph{grafted} simplicial complex in Section~\ref{grafting-section},
and show that for simplicial trees, being grafted and being unmixed
are equivalent conditions. The notion of grafting brings us to
Section~\ref{CM-section}, where we prove that grafted simplicial
complexes are Cohen-Macaulay, from which it follows that a simplicial
tree is unmixed if and only if it is Cohen-Macaulay.

\section{Definitions and notation}\label{basic-setup}
 
In this section we define the basic notions that we will use later in
the paper. Some of the proofs that appeared earlier in~\cite{Fa} have been
omitted here; we refer the reader to the relevant sections of~\cite{Fa}
when that is the case.

\begin{definition}[simplicial complex, facet and more] 
A \emph{simplicial complex} $\Delta$ over a set of vertices $V=\{ \vs
\}$ is a collection of subsets of $V$, with the property that $\{ v_i
\} \in \Delta$ for all $i$, and if $F \in \Delta$ then all subsets of
$F$ are also in $\Delta$ (including the empty set). An element of
$\Delta$ is called a \emph{face} of $\Delta$, and the \emph{dimension}
of a face $F$ of $\Delta$ is defined as $|F| -1$, where $|F|$ is the
number of vertices of $F$.  The faces of dimensions 0 and 1 are called
\emph{vertices} and \emph{edges}, respectively, and $\dimn \emptyset
=-1$.

The maximal faces of $\Delta$ under inclusion are called \emph{facets}
of $\Delta$. The dimension of the simplicial complex $\Delta$ is the
maximal dimension of its facets; in other words $$\dimn \Delta =\max
\{ \dimn F \st F \in \Delta \}.$$ 

We denote the simplicial complex $\Delta$ with facets $\Fs$ by
$$\Delta = \tuple{\Fs}$$ and we call $\{ \Fs \}$ the \emph{facet set}
of $\Delta$. 

A simplicial complex with only one facet is called a \emph{simplex}.

\end{definition}

\begin{definition}[subcollection] By a \emph{subcollection} of a simplicial
  complex $\Delta$ we mean a simplicial complex whose facet set is a
  subset of the facet set of $\Delta$.
\end{definition}

\begin{definition}[connected simplicial complex] A simplicial complex
  $\D=\tuple{\Fs}$ is \emph{connected} if for every pair $i,j$, $1
  \leq i < j \leq q$, there exists a sequence of facets 
  $$F_{t_1},\ldots,F_{t_r}$$
  of $\D$ such that $F_{t_1}=F_i$,
  $F_{t_r}=F_j$ and $$F_{t_s} \cap F_{t_{s+1}} \neq \emptyset$$
  for
  $s=1,\ldots,r-1$.

\end{definition}

An equivalent definition is stated on page~222 of~\cite{BH}: $\D$ as above
is \emph{disconnected} if its vertex set $V$ can be partitioned as
$V=V_1 \cup V_2$, where $V_1$ and $V_2$ are nonempty subsets of $V$,
such that no facet of $\D$ has vertices in both $V_1$ and $V_2$.
Otherwise $\D$ is \emph{connected}.

\idiot{ I call the sequence of facets in the definition above a
  ``path'' between $F_i$ and $F_j$.
  
  $\bullet$ To see that these two definitions are equivalent, label
  the first definition as (1) and the second one as (2).
  
  Let $V^i_1=\{V(F_j)|\ {\rm there \ is \ a \ path\ between\ } F_i
  {\rm \ and\ } F_j \}$, and let $V^i_2=V\setminus V^i_1$. Then obviously
  $V = V^i_1 \cup V^i_2$ and $V^i_1 \cap V^i_2 = \emptyset$.
  
  If $\D$ is disconnected in (1), then for some $i$ both $V^i_1$ and
  $V^i_2$ are nonempty, and so $\D$ is disconnected in (2).
  
  Suppose that $\D$ is disconnected in (2). So $V = V_1 \cup V_2$
  where both $V_1$ and $V_2$ are nonempty. Pick the facets $F_1$ and
  $F_2$ such that $V(F_1) \subseteq V_1$ and $V(F_2) \subseteq V_2$.
  If there is a path $H_0=F_1,H_1,\ldots,H_{t-1},H_t=F_2$ between
  $F_1$ and $F_2$, then by definition $H_i \cap H_{i+1} \neq
  \emptyset$, and so all the $H_i$ (including $F_1$ and $F_2$) have
  vertices in $V_1$ or all of them have vertices in $V_2$--
  contradiction. So there is no path between $F_1$ and $F_2$, which
  means that $\D$ is disconnected in (1).
  
  $\bullet$ The definition of ``path'' here is rather clumsy, I think.
  It does not replace the one in graph theory. For example, in the
  graph-tree $\tuple{xy,xz,xu}$, the sequence $xy,xz,xu$ is a path!
  This is not a ``walk''.}

\begin{definition}[facet ideal, non-face ideal] Let  $\Delta$
 be a simplicial complex over $n$ vertices labeled $\vs$. Let $k$ be a
field, $\xs$ be indeterminates, and $R$ be the polynomial ring
$k[\xs]$.

\noindent (a) We define $\F(\Delta)$ to be the ideal of $R$ generated
by all the square-free monomials $x_{i_1}\ldots x_{i_s}$, where
$\{v_{i_1},\ldots, v_{i_s}\}$ is a facet of $\Delta$. We call
$\F(\Delta)$ the \emph{facet ideal} of $\Delta$.

\noindent (b) We define $\N(\Delta)$ to be the ideal of $R$ generated
by all the square-free monomials $x_{i_1}\ldots x_{i_s}$, where
$\{v_{i_1},\ldots, v_{i_s}\}$ is not a face of $\Delta$. We call
$\N(\Delta)$ the \emph{non-face ideal} or the \emph{Stanley-Reisner
ideal} of $\Delta$.
\end{definition}

        We refer the reader to~\cite{S} and~\cite{BH} for an extensive
        coverage of the theory of Stanley-Reisner ideals.
        
        Throughout this paper we often use $\xs$ to denote both the
        vertices of $\D$ and the variables appearing in $\F(\D)$.

\begin{definition}[facet complex, non-face complex] Let
  $I=(\Ms)$ be an ideal in a polynomial ring $k[\xs]$, where $k$ is a
field and $\Ms$ are square-free monomials in $\xs$ that form a minimal
set of generators for $I$.

\noindent (a) We define $\df(I)$ to be the simplicial complex over a
set of vertices $\vs$ with facets $\Fs$, where for each $i$, $F_i=\{v_j
\st  x_j|M_i, \ 1 \leq j \leq n \}$. We call $\df(I)$ the \emph{facet
complex} of $I$.

\noindent (b) We define $\dn(I)$ to be the simplicial complex over a 
set of vertices $\vs$, where $\{v_{i_1},\ldots, v_{i_s}\}$ is a face
of $\dn(I)$ if and only if $x_{i_1}\ldots x_{i_s} \notin I$. We call
$\dn(I)$ the \emph{non-face complex} or the \emph{Stanley-Reisner
complex} of $I$.

\end{definition}

Facet ideals give a one-to-one correspondence between simplicial
complexes and square-free monomial ideals.

Notice that given a square-free monomial ideal $I$ in a polynomial
ring $k[\xs]$, the vertices of $\df(I)$ are those variables that
divide a monomial in the generating set of $I$; this set may not
necessarily include all elements of $\{\xs\}$. The fact that some
extra variables may appear in the polynomial ring does not affect the
algebraic or combinatorial structure of $\df(I)$. On the other hand,
if $\D$ is a simplicial complex, being able to consider the facet
ideals of its subcomplexes as ideals in the same ring simplifies many
of our discussions.

\begin{example}\label{example1} Let $\Delta$ be the simplicial complex below.

\[ \includegraphics{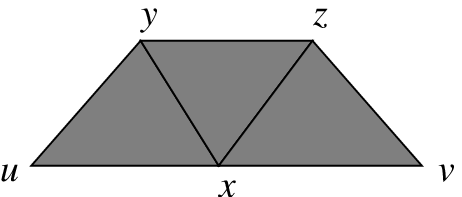} \]


Here $\N(\Delta)=(yv,zu,uv)$, $\F(\Delta)=(xyu,xyz,xzv)$ are
ideals in the polynomial ring $k[x,y,z,u,v]$. \end{example}

\begin{example}\label{example2} If $I=(xy,xz) \subseteq k[x,y,z]$, then 
$\dn(I)$ is the 1-dimensional simplicial complex:
\[ \includegraphics{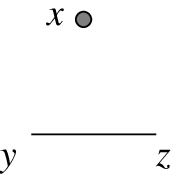}\] 
and $\df(I)$ is the simple graph
\[\includegraphics{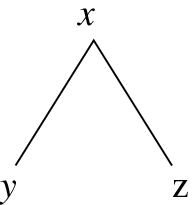}\] 
In this special case $I$ is also called the \emph{edge ideal} of the
graph $\df(I)$ (this terminology is due to Villarreal;
see~\cite{Vi1}). \end{example}

        We now generalize some notions from graph theory to simplicial
        complexes.

\begin{definition}[minimal vertex cover, vertex covering number, unmixed] 
  Let $\Delta$ be a simplicial complex with vertex set $V$ and facets
  $\Fs$. A \emph{vertex cover} for $\Delta$ is a subset $A$ of $V$,
  with the property that for every facet $F_i$ there is a vertex $v
  \in A$ such that $v \in F_i$. A \emph{minimal vertex cover} of
  $\Delta$ is a subset $A$ of $V$ such that $A$ is a vertex cover, and
  no proper subset of $A$ is a vertex cover for $\Delta$. The smallest
  cardinality of a vertex cover of $\Delta$ is called the \emph{vertex
    covering number} of $\Delta$ and is denoted by $\al(\D)$.

A simplicial complex $\Delta$ is \emph{unmixed} if all of its minimal
vertex covers have the same cardinality.
\end{definition}

Note that a simplicial complex may have several minimal vertex covers.

\begin{definition}[independent set, independence number]
  Let $\D$ be a simplicial complex. A set $\{F_1, \ldots, F_u\}$ of
  facets of $\D$ is called an \emph{independent set} if $F_i \cap F_j
  =\emptyset$ whenever $i \neq j$. The maximum possible cardinality of
  an independent set of facets in $\D$, denoted by $\beta(\D)$, is
  called the \emph{independence number} of $\D$. An independent set of
  facets which is not a proper subset of any other independent set is
  called a \emph{maximal independent set} of facets.
\end{definition}

\begin{example}\label{example11} If $\D$ is the simplicial complex 
\[ \includegraphics[height=1in]{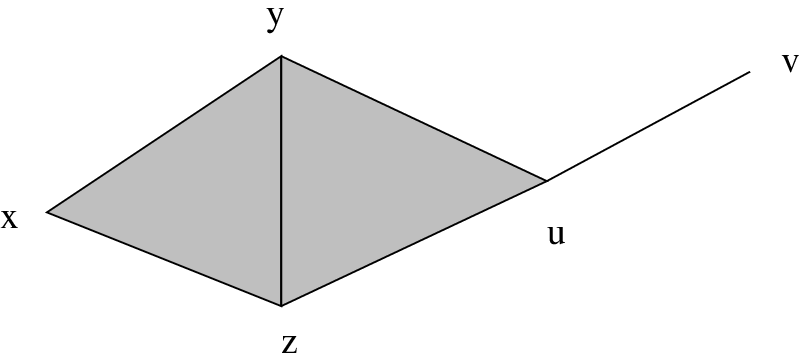} \]
then $\beta(\D)=2$. Also, $\D$ is unmixed as its minimal vertex
covers, listed below, all have cardinality equal to two:
$$\{x,u \}, \{y,u\}, \{y,v\}, \{z,u\}, \{z,v\}$$ 

This, by the way, is an example of a ``grafted'' tree (see
Definitions~\ref{tree} and~\ref{grafting}).  We show later in the
paper that all grafted trees are unmixed.

The graph $\df(I)$ in Example~\ref{example2} however is not unmixed.
This is because $\{x\}$ and $\{y,z\}$ are both minimal vertex covers
for $\df(I)$ of different cardinalities. In this case
$\al(\df(I))=\beta(\df(I))=1$. The same argument shows that the simplicial complex in Example~\ref{example1} is not unmixed. 

\end{example}

The following is an easy but very useful observation; see
Proposition~1 in~\cite{Fa} for a proof.

\begin{proposition}\label{minprime} Let $\Delta$ be a simplicial
  complex over $n$ vertices labeled $\xs$. Consider the ideal
  $I=\F(\Delta)$ in the polynomial ring $R=k[\xs]$ over a field $k$.
  Then an ideal $p=(x_{i_1},\ldots,x_{i_s})$ of $R$ is a minimal prime
  of $I$ if and only if $\{x_{i_1},\ldots,x_{i_s}\}$ is a minimal
  vertex cover for $\Delta$.
\end{proposition}

We say that a simplicial complex $\D$ over a set of vertices $\xs$ is
Cohen-Macaulay if for a given field $k$, the quotient ring
$$k[\xs]/\F(\D)$$ is Cohen-Macaulay. It follows directly from
Proposition~\ref{minprime}, or from an elementary duality with
Stanley-Reisner theory discussed in Corollary~2 of~\cite{Fa}, that in order
for $\D$ to be Cohen-Macaulay, it has to be unmixed.

\begin{proposition}[A Cohen-Macaulay simplicial complex is
 unmixed]\label{CM-is-unmixed} Suppose that $\Delta$ is a simplicial
complex with vertex set $\xs$. If $k[\xs]/\F(\Delta)$ is
Cohen-Macaulay, then $\Delta$ is unmixed.\end{proposition}

\begin{discussion}\label{dimension} It is worth observing that for a 
square-free monomial ideal $I$, there is a natural way to construct
$\dn(I)$ and $\df(I)$ from each other using the structure of the
minimal primes of $I$. To do this, consider the vertex set $V$
consisting of all variables that divide a monomial in the generating
set of $I$. The following correspondence holds:

        \begin{center} $F=$ facet of $\dn(I)$ $\longleftrightarrow $
        $V \setminus F=$ minimal vertex cover of $\df(I)$ \end{center}

Also $$I =\bigcap p$$ where the intersection is taken over all
prime ideals $p$ of $k[V]$ that are generated by a minimal vertex
cover of $\df(I)$ (or equivalently, primes $p$ that are generated by
$V \setminus F$, where $F$ is a facet of $\dn(I)$; see~\cite{BH} Theorem 
5.1.4). 

        Regarding the dimension and codimension of $I$, note that by
        Theorem 5.1.4 of~\cite{BH} and the discussion above, setting
        $R=k[V]$ as above, we have
        $$\dimn R/I =\dimn \dn(I) +1 = |V| -{\rm \ vertex\ covering\
number\ of\ }\df(I)$$ and $$\height I\ ={\rm \ vertex\ covering\
number\ of\ } \df(I).$$

\idiot{ This is because $\dimn \dn(I) = \dimn F = |F|-1$, where $F$ is
the highest dimensional facet of $\dn(I)$. So $\dimn R/I = |F| = |V|
-{\rm \ minimal\ vertex\ covering\ number\ of\ }\df(I)$.}

\end{discussion}

We illustrate all this through an example.

\begin{example} For $I=(xy, xz)$, where $\df(I)$ and $\dn(I)$ are drawn 
  in Example~\ref{example2}, we have:

\begin{center}
\begin{tabular}{cc}
\hspace{.3in}\underline{facets of $\dn(I)$}\hspace{.3in} &
\hspace{.3in}\underline{minimal vertex covers of
$\df(I)$}\hspace{.3in} \\ &\\ $\{ x \}$ & $\{ y, z \}$ \\ $\{ y, z
\}$& $\{ x \}$
\end{tabular}
\end{center}
 Note that $I = (x) \cap (y,z)$, and $$\dimn k[x,y,z]/(xy,xz)=2$$
 as  asserted in Discussion~\ref{dimension} above.

\end{example}

       A notion crucial to the rest of the paper is ``removing a
        facet''. We want the removal of a facet from a simplicial
        complex to correspond to dropping a generator from its facet
        ideal. We record this definition.

\begin{definition}[facet removal]\label{removal} Suppose $\Delta$
  is a simplicial complex with facets $\Fs$ and $\F(\Delta)=(\Ms)$ its
  facet ideal in $R=k[\xs]$. The simplicial complex obtained by
  \emph{removing the facet} $F_i$ from $\Delta$ is the simplicial
  complex
 $$\D \rmv{F_i}=\tuple{F_1,\ldots,\hat{F}_{i},\ldots,F_q}.$$
 \end{definition}
 
 Note that $\F(\D \rmv{F_i}) = (M_1,\ldots,\hat{M}_{i} ,\ldots,M_q)$.
 
 Also note that the vertex set of $\D \rmv{F_i}$ is a subset of the
 vertex set of $\D$.

\begin{example} let $\D$ be the simplicial complex in
 Example~\ref{example11} with facets $F=\{x,y,z\}$, $G=\{y,z,u\}$ and
 $H=\{u,v\}$. Then $\D \rmv{F} =\tuple{G,H}$  is a simplicial complex
 with vertex set $\{y,z,u,v\}$. 
\end{example}


\section{Trees}\label{tree-section}

In~\cite{Fa} we extended the notion of a ``tree'' from graphs to simplicial
complexes. The construction, at the time, was motivated by two
factors: the restriction to graphs should produce the classic
graph-theoretical definition of a tree, and the new structure should
fit into a machinery that proves that the facet ideal of a tree
satisfies Sliding Depth condition (Theorem~1 of~\cite{Fa}).

The resulting definition not only satisfies those two properties, but
as we prove later in this paper, it also generalizes graph-trees in
the sense of Cohen-Macaulayness, which confirms that algebraically
this in fact is the optimal way to extend the definition of a tree.

         Recall that a connected graph is a tree if it has no cycles;
         for example, a triangle is not a tree. An equivalent
         definition states that a connected graph is a tree if every
         subgraph has a \emph{leaf}, where a leaf is a vertex that
         belongs to only one edge of the graph. This latter
         description is the one that we adapt, with a slight change in
         the definition of a leaf, to the class of simplicial
         complexes.

\begin{definition}[leaf, joint, universal set]~\label{leaf} 
Suppose that $\Delta$ is a simplicial complex. A facet $F$ of $\Delta$
 is called a \emph{leaf} if either $F$ is the only facet of $\Delta$,
 or there exists a facet $G$ in $\Delta\rmv{F}$, such that $$ F \cap
 F' \subseteq F \cap G$$ for every facet $F' \in \Delta\rmv{F}$.
 
In other words, $F$ is a leaf of $\D$ if it intersects $\D \rmv{F}$ in
a face of $\D \rmv{F}$.

 The set of all $G$ as above is denoted by $\U_\Delta (F)$ and called
 the \emph{universal set} of $F$ in $\Delta$. If $G \in \U_\Delta (F)$
 and $F\cap G \neq \emptyset$, then $G$ is called a \emph{joint} of
 $F$.

\end{definition}

Another way to describe a leaf is the following: (with assumptions as
above) $F$ is a leaf if either $F$ is the only facet of $\D$ or the
intersection of $F$ with the simplicial complex $\D \rmv{F}$ is a face
of $\D \rmv{F}$.

\begin{definition}[free vertex]\label{free vertex} A vertex of a simplicial 
complex $\D$ is \emph{free} if it belongs to exactly one facet of
$\D$. \end{definition}

In order to be able to quickly identify a leaf in a simplicial
complex, it is important to notice that a leaf must have a free
vertex. This follows easily from Definition~\ref{leaf}: otherwise, a
leaf $F$ would be contained in its joints, which would contradict the
fact that a leaf is a facet.

\begin{example}\label{non-leaf-example} 

The simplicial complex in example~\ref{example1} has two leaves:
$\{x,y,u\}$ and $\{x,z,v\}$. The one below has no leaves,
because every vertex is shared by at least two facets.

\[ \includegraphics[height=.6in]{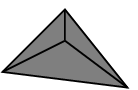} \]
\end{example}

\begin{example}\label{good-example} In the simplicial complex below 
with facets $F_1=\{a,b,c \}$, $F_2=\{ a,c,d \}$ and $F_3=\{ b,c,d,e
  \}$, the only candidate for a leaf is the facet $F_3$ (as it is the
  only facet with a free vertex), but neither one of $F_1 \cap F_3$ or
  $F_2 \cap F_3$ is contained in the other, so there are no leaves.

\[ \includegraphics{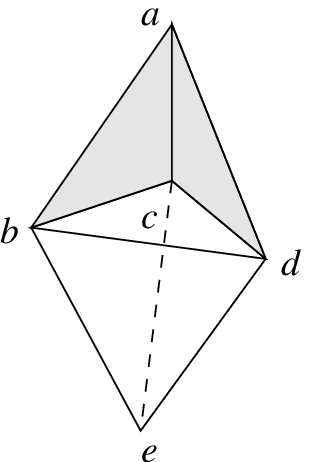} \]

\end{example}

\begin{definition}[tree]~\label{tree} Suppose that $\Delta$ 
  is a connected simplicial complex.  We say that $\Delta$ is a
  \emph{tree} if every nonempty subcollection of $\Delta$ (including
  $\Delta$ itself) has a leaf.
  
  Equivalently, a connected simplicial complex $\Delta$ is a tree if
  every nonempty \emph{connected} subcollection of $\Delta$ has a
  leaf.
 \end{definition}

\begin{definition}[forest]\label{forest} A simplicial complex
  $\Delta$ with the property that every connected component of
  $\Delta$ is a tree is called a \emph{forest}. In other words, a
  forest is a simplicial complex with the property that every nonempty
  subcollection has a leaf.
\end{definition}

The simplicial complex in Example~\ref{example1} above is a tree,
whereas the ones in Examples~\ref{non-leaf-example}
and~\ref{good-example} are not, as they have  no leaves.

Here is a slightly less straightforward example:

\begin{example}

The simplicial complex on the left is not a tree, because although all
three facets $\{x,y,u\}$, $\{x,v,z\}$ and $\{y,z,w\}$ are leaves, if
one removes the facet $\{x,y,z\}$, the remaining simplicial complex
(on the right) has no leaf.

\begin{center}
\begin{tabular}{llr}
&&\\
  \parbox{1in}{\epsfig{file=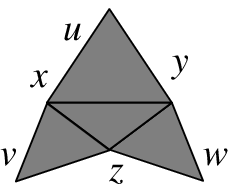, height=.8in}} 
& $\stackrel{\mbox{remove $\{x,y,z\}$}}{\vector(1,0){100}}$ 
& \parbox{1in}{\epsfig{file=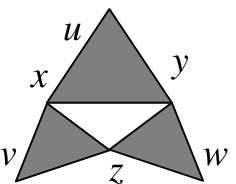, height=.8in}}\\
&&
\end{tabular}
\end{center}

\end{example}

        Notice that in the case that $\Delta$ is a graph,
        Definition~\ref{tree} agrees with the definition of a tree in
        graph theory, with the difference that now the term ``leaf''
        refers to an edge, rather than a vertex.


\section{Basic properties of trees}\label{basic-trees-section}

\begin{lemma}[A tree has at least two leaves]\label{two leaves}
 Let $\Delta$ be a tree of two or more facets. Then $\Delta$ has at
 least two leaves. \end{lemma}

        \begin{proof} Suppose that $\Delta$ has $q$ facets $\Fs$ where
        $q \geq 2$. We prove the lemma by induction on $q$.
        
        The case $q=2$ follows from the definition of a leaf.

        To prove the general case suppose that $F_1$ is a leaf of
        $\Delta$ and $G_1 \in \U_\Delta (F_1)$.  Consider the
        subcomplex $\Delta '= \tuple{F_2, \ldots, F_q}$ of
        $\Delta$. By induction hypothesis $\Delta '$ has two distinct
        leaves; say $F_2$ and $F_3$ are those leaves. At least one of
        $F_2$ and $F_3$ must be different from $G_1$; say $F_2 \neq
        G_1$.  We show that $F_2$ is a leaf for $\Delta$.

        Let $G_2 \in \U_{\Delta '} (F_2)$. Given any
        facet $F_i$ with $i \neq 1, 2$, we already know by the fact
        that $F_2$ is a leaf of $\Delta '$
        $$ F_i \cap F_2 \subseteq G_2 \cap F_2 .$$ We need to verify
        this for $i=1$.

         Since $F_1$ is a leaf for $\Delta$ and $F_2 \neq F_1$,
        $$ F_2 \cap F_1 \subseteq G_1 \cap F_1 .$$ Intersecting both
        sides of this inclusion with $F_2$, we obtain
        $$ F_2 \cap F_1 \subseteq G_1 \cap F_1 \cap F_2 \subseteq G_1
        \cap F_2 \subseteq G_2 \cap F_2$$
        where the last inclusion holds
        because $G_1 \neq F_2$ and $F_2$ is a leaf of
        $\Delta '$.

        It follows that $F_2$, as well as $F_1$, is a leaf for $\Delta$.

        \end{proof}

        A promising property of trees from an algebraic point of view
        is that they behave well under localization, i.e. the
        localization of a tree is a forest. This property is in
        particular useful when making inductive arguments on trees, as
        localization usually corresponds to reducing the number of
        vertices of a simplicial complex.  Before proving this, we
        first determine what the localization of a simplicial complex
        precisely looks like.

\begin{discussion}[On the localization of a simplicial
  complex]\label{local-discussion} Suppose that 
  $$\D=\tuple{\Fs}$$
  is a simplicial complex over the vertex set
  $V=\{\xs\}$. Let $p$ be a prime ideal of $k[\xs]$ generated by a
  subset of $\{\xs\}$ that contains $I=\F(\D)$ (We show later in the
  proof of Lemma~\ref{localization} that this is the main case that we
  need to study). We would like to see what the simplicial complex
  associated to $I_p$ looks like. 
  
  So $$p=(x_{i_1},\ldots, x_{i_r}).$$
  Now suppose
  $$I=(\Ms)$$
  where each $M_i$ is the monomial corresponding to the
  facet $F_i$. It follows that
  $$I_p=(M_1', \ldots, M_q')$$
  where each $M_i'$ is obtained by
  dividing $M_i$ by the product of all the variables in $V
  \setminus\{x_{i_1},\ldots, x_{i_r}\}$ that appear in $M_i$. Some of
  the monomials in the generating set of $I_p$ are redundant after
  this elimination, so without loss of generality we can write:
\begin{eqnarray}\label{e:local} I_p=(M_1', \ldots, M_t') \end{eqnarray}
 where $M'_{t+1},\ldots, M'_q$ are the redundant monomials.

We use the notation $\df(I_p)$ to indicate the simplicial complex
  associated to the monomial ideal with the same generating set as the
  one described for $I_p$ in (\ref{e:local}),  in the polynomial
  ring $k[x_{i_1},\ldots, x_{i_r}]$. It follows that:
$$\df(I_p)=\tuple{F'_1, \ldots, F'_t}$$ where for each $i$, $$F'_i=F_i
\cap \{ x_{i_1},\ldots, x_{i_r}\}$$ and $F'_{t+1}, \ldots, F'_q$ each
contain at least one of $F'_1, \ldots, F'_t$. This simplicial complex
is called the \emph{localization} of $\D$ at the prime ideal $p$.

Note that every minimal vertex cover $A$ of $\D$ that is contained in
$\{ x_{i_1},\ldots, x_{i_r}\}$ remains a minimal vertex cover of
$\df(I_p)$, as the minimal prime over $I$ generated by the elements of
$A$ remains a minimal prime of $I_p$. 

Moreover if $\D$ is unmixed then $\df(I_p)$ is also unmixed.
Algebraically, this is easy to see, as the height of the minimal
primes of $I_p$ remain the same. One can also see it from a
combinatorial argument: If $B \subseteq \{x_{i_1},\ldots, x_{i_r}\}$
is a minimal vertex cover for $\df(I_p)$, then $B$ covers all facets
$F'_1, \ldots, F'_t$, and therefore $F'_{t+1}, \ldots, F'_q$, as well.
Therefore $B$ covers all of $F_1, \ldots, F_q$, and so has a subset
$B'$ of cardinality $\al(\D)$ that is a minimal vertex cover for $\D$,
and so $B'$ must cover $\df(I_p)$ as well. Therefore $B'=B$.

We have thus shown that:

\begin{lemma}[Localization of an unmixed simplicial complex is unmixed]\label{unmixed-localizes} Let $\Delta$ be an unmixed simplicial complex with vertices 
  $\xs$, and let $I= \F (\Delta)$ be the facet ideal of $\Delta$ in
  the polynomial ring $R=k[\xs]$ where $k$ is a field. Then for any
  prime ideal $p$ of $R$, $\df (I_p)$ is unmixed with $\al(\df
  (I_p))=\al(\D)$. \end{lemma}

\end{discussion}

We examine a specific case:

\begin{example}\label{local-example} Let $\D$ be the simplicial complex below
  with $I=(xyu,xyz,xzv)$ its facet ideal in the polynomial ring
  $R=k[x,y,z,u,v]$.

\[ \includegraphics{labeled-tree-or-not.eps} \]

Let $p=(u,x,z)$ be a prime ideal of $R$. Then $I_p=(xu,xz,
xz)=(xu,xz)$. The tree $\df(I_p)$, shown below, has minimal vertex
covers $\{ x\}$ and $\{u,z\}$, which are the generating sets for the
minimal primes of $I_p$.

\[ \includegraphics{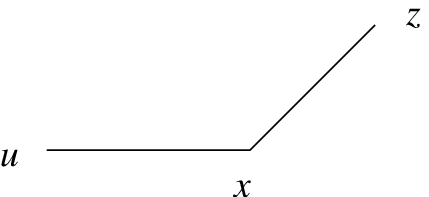} \]

If $q= (y,z,v)$ then $I_q=(y,yz,zv)=(y,zv)$ which corresponds to the
forest $\df(I_q)$ drawn below with minimal vertex covers $\{y,z\}$ and
$\{y,v\}$.

\[ \includegraphics{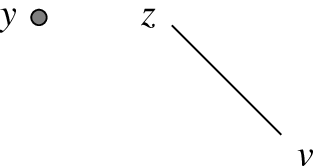} \]
\end{example}

Example~\ref{local-example} above also demonstrates the following
lemma.

\idiot{ The reason I left the proof of this lemma here in this paper
  is that in the previous paper I had the wrong idea about
  ``redundancy'', and so I describe it more carefully here}

\begin{lemma}[Localization of a tree is a forest]\label{localization}
  Let $\Delta$ be a tree with vertices $\xs$, and let $I= \F (\Delta)$
  be the facet ideal of $\Delta$ in the polynomial ring $R=k[\xs]$
  where $k$ is a field. Then for any prime ideal $p$ of $R$, $\df
  (I_p)$ is a forest. \end{lemma}

        \begin{proof} The first step is to show that it is enough to
        prove this for prime ideals of $R$ generated by a subset of
        $\{ \xs \}$. To see this, assume that $p$ is a prime ideal of
        $R$ and that $p'$ is another prime of $R$ generated by all
        $x_i \in \{ \xs \}$ such that $x_i \in p$ (recall that the
        minimal primes of $I$ are generated by subsets of $\{ \xs
        \}$). So $p' \subseteq p$. If $I=(\Ms)$, then $$I_{p'} =
        ({M_1}',\ldots,{M_q}')$$ where for each $i$, ${M_i}'$ is the
        image of $M_i$ in $I_{p'}$. In other words, ${M_i}'$ is
        obtained by dividing $M_i$ by the product of all the $x_j$
        such that $x_j| M_i$ and $x_j \notin p'$. But $x_j \notin p'$
        implies that $x_j \notin p$, and so it follows that ${M_i}'
        \in I_p$. Therefore $I_{p'} \subseteq I_p$. On the other hand
        since $p' \subseteq p$, $I_p \subseteq I_{p'}$, which implies
        that $I_{p'} = I_p$ (the equality and inclusions of the ideals
        here mean equality and inclusion of their generating sets).

        We now prove the theorem for $p=(x_{i_1},\ldots, x_{i_r})$.
        Following the setup in Discussion~\ref{local-discussion}, we
        let $$\Delta= \tuple{\Fs}$$
        $$I_p=(M_1', \ldots, M_t')$$
        $$\D'=\df(I_p)=\tuple{F'_1, \ldots, F'_t}.$$
        for some $t \leq
        q$.  
        
        To show that $\D'$ is a forest, we need to show that every
        nonempty subcollection of $\D'$ has a leaf.
        
        Let $$\Delta_1 '=\tuple{ F_{j_1} ',\ldots, F_{j_s} '}$$ be a
        subcollection of $\Delta '$ where $F_{j_1} ',\ldots, F_{j_s}
        '$ are distinct facets. If $s=1$, $F'_{j_1}$ is obviously a
        leaf and so we are done; so suppose $s>1$. Consider the
        corresponding subcollection
        $$\Delta_1=\tuple{ F_{j_1} ,\ldots, F_{j_s} }$$ of $\Delta$,
        which has a leaf, say $F_{j_1}$. So there exists $G \in
        \Delta_1 \rmv{F_{j_1}}$, such that $$ F_{j_1} \cap F \subseteq
        F_{j_1} \cap G$$ for every facet $F \in \tuple{F_{j_2}
        ,\ldots, F_{j_s} }$. Now since each of the $ F_{j_u}'$ is a
        nonempty facet of $\Delta_1 '$ and $G' \neq F'_{j_1}$, the
        same statement holds in $\Delta_1 '$; so
        $$ F_{j_1}' \cap F' \subseteq F_{j_1}' \cap G'$$ for every
        facet $F' \in \Delta_1' \rmv{ F_{j_1}'}$.
        This implies that $ F_{j_1}'$ is a leaf for $\Delta_1 '$.

        \end{proof}

      
\section{Simplicial complexes as higher dimensional 
graphs}\label{n-partite-section}

In this section we study simplicial complexes as graphs with higher
dimension, drawing results that will help us later in inductive
arguments on unmixed trees. 

\void{We also include, independently of the rest
of the paper, a brief discussion on multi-partite simplicial
complexes. This is inspired by graph theory, where K\"onig's theorem
(generalized in Theorem~\ref{konig} below) holds for all bipartite
graphs. The natural extension of ``bipartite'' to higher dimension,
however, does not provide suitable conditions for Theorem~\ref{konig}
to hold. On the other hand, we can show that all trees are
multi-partite (Theorem~\ref{theorem-tree-multipartite}), which gives
us much insight into the vertex structure of trees.
}

\begin{lemma}\label{unmixed-induction} If $\D$ is a simplicial complex  
  that has a leaf $F$ with joint $G$, then $\al(\D\rmv{G})=\al(\D)$.
\end{lemma}

        \begin{proof} Suppose $\al(\D)=r$. Let $\D' = \D\rmv{G}$ and
          let $A$ be a vertex cover of minimal cardinality for $\D'$,
          which implies that $|A| \leq r$, as any vertex cover of $\D$
          has a subset that is a vertex cover of $\D'$. Since $F$ is
          a facet of $\D'$, there exists a vertex $x \in A$ that
          belongs to $F$.  If $x$ is a free vertex of $F$, we may
          replace it by a non-free vertex of $F$ to get a vertex cover
          $A''$ of $\D'$, with a subset $A'$ that is a minimal vertex
          cover of $\D'$, and so $|A'| \leq |A|$.  But now $A'$ is a
          minimal vertex cover for all of $\D$, and so $|A'|=|A|=r$
          which implies that $\al(\D')= \al(\D)=r$.\end{proof}

\begin{corollary}\label{unmixed-corollary} If the simplicial complex
$\D$ is a tree and $G \in \D$ is a joint, then $\al(\D\rmv{G})=
  \al(\D)$. \end{corollary}

This means that in a tree with more than one facet, it is always
possible to remove a facet without reducing the vertex covering
number. Moreover we show in Proposition~\ref{joint-removal-unmixed}
that if $\D$ is an unmixed tree with a joint $G$, then $\D \rmv{G}$ is
also unmixed. As a result, one can use induction on the number of
facets of an unmixed tree. Note that all these arguments remain valid
for a forest.

We are now ready to extend K\"{o}nig's theorem from graph theory.

\begin{theorem}[A generalization of K\"{o}nig's theorem]\label{konig} 
  If $\D$ is a simplicial complex that is a tree (forest) and
  $\al(\D)=r$, then $\D$ has $r$ independent facets, and therefore
  $\al(\D)=\be(\D)=r$.
\end{theorem}

        \begin{proof} We use induction on the number of facets $q$ of $\D$.
        If $q=1$, then there is nothing to prove since $\al(\D)=\be(\D)=1$.
        
        Suppose that the theorem holds for forests with less than $q$
        facets and let $\D$ be a forest with $q$ facets. If every
        connected component of $\D$ has only one facet, our claim
        follows immediately. Otherwise, by
        Corollary~\ref{unmixed-corollary} one can remove a joint of
        $\D$ to get a forest $\D'$ with $\al(\D') =r$, and so by
        induction hypothesis $\D'$ has $r$ independent facets, which
        are also independent facets of $\D$; so $\al(\D) \leq
        \be(\D)$. On the other hand it is clear that $\al(\D) \geq
        \be(\D)$, and so the assertion follows.
      \end{proof}

\section{The structure of an unmixed tree}\label{structure-theorem-section}
 
This section is the combinatorial core of the paper. Here we give a
precise description of the structure of an unmixed tree. It turns out
that a tree is unmixed if and only if it is ``grafted'' (see
Definition~\ref{grafting}). The notion of grafting is what
eventually builds a bridge between unmixed and Cohen-Macaulay
trees.

Below $V(\D)$ stands for the vertex set of $\D$.

\begin{lemma}\label{basic-lemma} Let $\D$ be an unmixed simplicial complex. 
  Suppose that $\al(\D) = \beta(\D)=r$, and $\{F_1, \ldots, F_r\}$ is
  a maximal independent set of facets of $\D$. Then every vertex of
  $\D$ belongs to one of the $F_i$. In other words, the vertex set of
  $\D$ is the disjoint union of the vertex sets of the $F_i$:
  $$V(\D)=V(F_1)\cup \ldots \cup V(F_r)$$
\end{lemma}

        \begin{proof} Let $x$ be an vertex of $\D$ that does not belong
          to any of the $F_i$. Then one can find a minimal vertex
          cover $A$ of $\D$ containing $x$ (this is always possible).
          But then $A$ must contain one vertex of each of the $F_i$ as
          well, which implies that $|A| \geq r+1$. Since $\D$ is
          unmixed, this is not possible. \end{proof}
        
        \idiot{$\bullet$ If $x$ is a vertex of $\D=\tuple{\Fs}$, it is
          always possible to build a minimal vertex cover of $\D$ that
          contains $x$. Suppose that $x \in F_1$, and let
          $F_i'=F_i\setminus F_1$. Now take a vertex cover consisting
          of vertices of $F_2', \ldots, F_q'$ and $x$. Any minimal
          vertex cover of $\D$ contained in this set has to contain
          $x$.}

\begin{remark}  Lemma~\ref{basic-lemma} does not hold in
          general for any unmixed simplicial complex. Take, for
          example, the case of a complete graph $G$ over 5 vertices
          labeled $x,y,z,u,v$ (every pair of vertices of $G$ are
          connected by an edge).  This graph is unmixed with
          $\al(G)=4$ and $\beta(G)=2$. However, $\{xy,uv\}$ is a
          maximal independent set of facets and the fifth vertex $z$
          of $G$ is missing from the vertex set of the graph
          $\tuple{xy,uv}$, which contradicts the claim of
          Lemma~\ref{basic-lemma}.\end{remark}

Lemma~\ref{basic-lemma} along with Theorem~\ref{konig} provides us
with the following property for unmixed trees.

\begin{corollary}\label{basic} If $\D$ is an unmixed tree with 
  $\al(\D)=r$, and $\{F_1, \ldots, F_r\}$ is a maximal
  independent set of facets of $\D$, then $V(\D)=V(F_1)\cup \ldots
  \cup V(F_r)$.
\end{corollary}

\begin{corollary}\label{leaf-independence} If $\D$ is an unmixed 
  tree, then any maximal independent set of facets of
  cardinality $\al(\D)$ of $\D$ contains all the leaves.  In
  particular, the leaves of an unmixed tree are independent.
\end{corollary}

        \begin{proof} Every leaf has a free vertex, and so it follows
          from above that a independent set of facets of cardinality
          $\al(\D)$ must contain all the leaves. The claim then
          follows. \end{proof}

\begin{corollary}\label{joint-not-leaf} If $\D$ is an unmixed 
  tree, then a maximal independent set of facets of
  cardinality $\al(\D)$ of $\D$ cannot contain a joint. In
  particular, a joint of an unmixed tree cannot be a leaf.
\end{corollary}

        \begin{proof} If $G$ is a joint, it has to intersect a leaf
          $F$ by definition, and as $F$ is in every
          maximal independent set of facets of cardinality
          $\al(\D)$, $G$ cannot be in any.
        \end{proof}

        But even more is true. For an unmixed tree $\D$, there is only
        one maximal independent set of facets with $\al(\D)$ elements,
        and that is the set consisting of all the leaves.  We prove
        this in Theorem~\ref{structuretheorem}.

The proposition below allows us to use induction on the number of
facets of an unmixed tree.

\begin{proposition}\label{joint-removal-unmixed}  Let $\D$ be an unmixed 
  tree with a leaf $F$, and let $G$ be a joint of $F$. Then $\D' = \D
  \rmv{G}$ is also unmixed.
\end{proposition}

          \begin{proof} We use induction on the number of vertices of $\D$.
            Let $$\D =\tuple{\Fs}$$ and $$V=\{ \xs \}$$ be the vertex
            set for $\D$.

            The case $n=1$ is clear.
 
            Suppose that $\al(\D)=r$ and $A$ is a minimal vertex cover
            for $\D'$.  By Corollary~\ref{unmixed-corollary}
            $\al(\D')=r$ as well.  If $A$ contains any vertex of $G$,
            then it is also a minimal vertex cover for $\D$ and hence
            of cardinality $r$. So suppose that $$A \cap G=\emptyset \ 
            \ {\rm and} \ \ |A|>r.$$

         {\bf Claim:} \emph{There is a vertex $x \in V \setminus (A
            \cup G)$.}
        \vspace{.2in}

        \begin{itemize}

        \item[ ]\noindent \emph{Proof of Claim:} If not, then
            \begin{eqnarray}\label{A=V-G} V=A \cup G. 
            \end{eqnarray}

            We show that this is not possible.
            
            Notice that for any $y \in A$ there is a facet $H \in \D'$
            such that $H \cap A=\{y\}$ (If no such $H$ existed, then
            $A \setminus \{y\}$ would also be a vertex cover).

            From~(\ref{A=V-G}) it follows that
             \begin{eqnarray}\label{HGH}
             H=(G \cap H) \cup \{y\}. \end{eqnarray}
            
            On the other hand using Theorem~\ref{konig} we can assume
            $\{F_1, \ldots, F_r\}$ is a maximal independent set of
            facets in $\D$.  By Corollary~\ref{joint-not-leaf} $$G
            \notin \{F_1, \ldots, F_r\}.$$
            As $|A| >r$, one of the
            $F_i$, say $F_r$, has to contain more than one element of
            $A$, so suppose
            $$A \cap F_r=\{y_1,\ldots,y_s\}$$ where $s>1$ and
            $y_1,\ldots,y_s$ are distinct elements of $A$.  It follows
            from~(\ref{A=V-G}) that \begin{eqnarray}\label{FrG} F_r =
            (F_r\cap G)\cup \{y_1,\ldots,y_s\}.  \end{eqnarray}
            
            From the discussion preceding (\ref{HGH}) above, one can
            pick $H_1, \ldots,H_s$ to be facets of $\D'$ such that
            \begin{eqnarray}\label{HiG} H_i=(G \cap H_i) \cup \{y_i\} 
            \end{eqnarray} for $i=1,\ldots,s$, and
            consider the tree
            $$\tuple{G, F_r,H_1, \ldots,H_s}$$ which by Lemma~\ref{two
            leaves} is supposed to have two leaves. But based on the
            descriptions of $F_r,H_1, \ldots,H_s$ in (\ref{FrG}) and
            (\ref{HiG}), only one facet of this tree, namely $G$,
            could possibly have a free vertex, which is a
            contradiction. This proves the claim.
           \end{itemize}

            We now proceed to showing that $|A|>r$ is not possible.
           
            Let $x \in V \setminus (A \cup G).$ We localize at the
            prime ideal $p$ generated by $V\setminus \{x\}$, and use
            the induction hypothesis.

            Let $$I =\F(\D) \ \ {\rm and} \ \ I'=\F(\D')$$ and let
            $$\D_p=\df(I_p) \ \ {\rm and} \ \ {\D'}_p=\df({I'}_p).$$ 
            
            From Discussion~\ref{local-discussion} we know that,
            without loss of generality, for some $t \leq q$
            $$\D_p= \tuple{\tilde{F}_1, \ldots, \tilde{F}_t}$$ where
            $\tilde{F}_i = F_i \setminus \{x\}$, and each of
            $\tilde{F}_{t+1}, \ldots, \tilde{F}_q$ contains at least
            one of $\tilde{F}_1, \ldots, \tilde{F}_t$.

            We also know by Lemma~\ref{localization} that $\D_p$
            is a forest whose vertex set is a proper subset of $V$. 
            
            By Lemma~\ref{unmixed-localizes} $\D_p$ is unmixed
            with $\al(\D_p)=r$ 
            
            We now focus on ${\D'}_p$. Besides possibly $\tilde{G}$,
            all other facets of $\D_p$ and ${\D'}_p$ are the same.  We
            show why this is true.

            Let $\tilde{F}_i \in {\D'}_p$. Then clearly $$ \tilde{F}_j
            \not\subseteq \tilde{F}_i {\rm \ for\ all\ } F_j \in \D',\
            j \neq i.$$ On the other hand, as $\tilde{G}=G$ and $G
            \not \subseteq F_i$, we have
            $$\tilde{G} \not\subseteq \tilde{F}_i$$ and so 
            $\tilde{F}_i \in \D_p.$
            
            Conversely, if $\tilde{F}_i \in \D_p$, then 
            $$\tilde{F}_j \not\subseteq \tilde{F}_i {\rm \ for\ all\ }
            F_j \in \D, \ j \neq i$$ which implies the same for all
            $F_j \in \D'$, and therefore $\tilde{F}_i \in {\D'}_p$.

            So there are two possible scenarios:

            \begin{enumerate}

              \item[\emph{Case 1.}] If $\tilde{G} \notin \D_p$, then
              $$\D_p={\D'}_p$$ which implies that $A$ is also a
              minimal vertex cover of $\D_p$, which is unmixed, and hence 
              $|A|=r$; a contradiction.

              \item[\emph{Case 2.}]  If $\tilde{G} \in \D_p$ then
              $$\tilde{F} \in \D_p.$$ If not, then for some facet $H$
              of $\D$, we have $\tilde{H} \subseteq \tilde{F}$, so $H
              \cap F \neq \emptyset$ and therefore, since $G$ is a
              joint of the leaf $F$,  $$H \cap F \subseteq G \cap F$$
              which immediately results in
              $$\tilde{H} \subseteq \tilde{G}$$ which is not possible.
              
              In fact, $\tilde{F}$ remains a leaf in $\D_p$, since if
              $\tilde{H}$ is a facet of $\D_p$ such that $\tilde{H}
              \cap \tilde{F} \neq \emptyset$, then $$\emptyset \neq
              H\cap F \subseteq G\cap F \Longrightarrow \tilde{H} \cap
              \tilde{F} \subseteq \tilde{G} \cap \tilde{F} $$ and so
              $\tilde{G}$ is a joint of $\D_p$.

              Now by the induction hypothesis, $${\D'}_p=\D_p
              \rmv{\tilde{G}}$$ is an unmixed forest. This again implies that
              $|A|=r$; a contradiction.
          \end{enumerate}
\end{proof}

\begin{example} Although  not obvious at a first glance, 
  Proposition~\ref{joint-removal-unmixed} does not necessarily hold if
  $G$ is not a tree. The following example of an unmixed graph $G$
  with a leaf demonstrates this point.

\[ \includegraphics{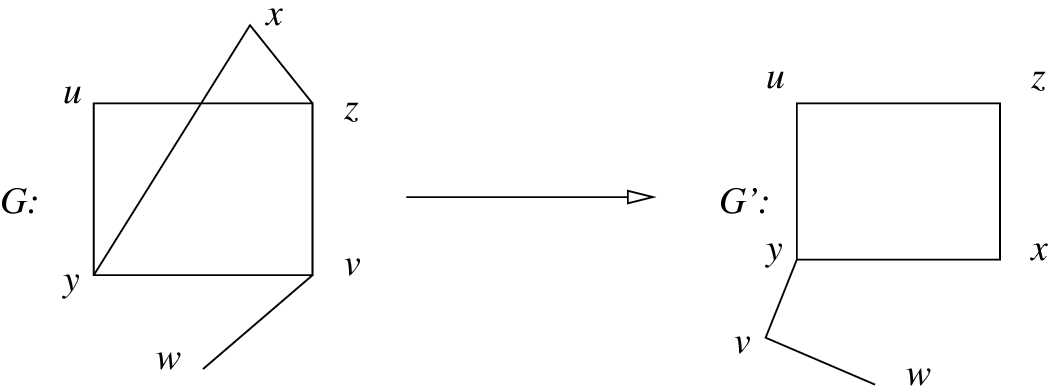} \]

The graph $G$ above was taken from the table of unmixed graphs
in~\cite{Vi2}. The minimal vertex covers of $G$, all of cardinality 3,
are $\{w,z,y\}, \{v,x,u\},$ and $ \{v, z, y\}$. But once one removes
the joint $\{v,z\}$, $G'$ has minimal vertex covers $ \{w,y,z\}$ and
$\{w,y,x,u\}$ of different cardinalities, and is therefore not
unmixed.
\end{example}

\begin{theorem}[Structure theorem for unmixed trees]\label{structuretheorem}
  Suppose that $\D$ is an unmixed tree with more than one facet such
  that $\al(\D)=r$. Then $\D$ can be written as
$$\D= \tuple{F_1, \ldots, F_r} \cup \tuple{G_1, \ldots, G_s}$$ with
the following properties:

\begin{enumerate}
\item[(i)] $F_1, \ldots, F_r$ are all the leaves of $\D$;
\item[(ii)] $\{G_1, \ldots, G_s\} \cap \{F_1, \ldots, F_r\} = \emptyset$;
\item[(iii)] For $i \neq j$, $F_i \cap F_j = \emptyset$;
\item[(iv)] If a facet $H \in \D$ is not a leaf, then it does not
contain a free vertex.
\end{enumerate}
\end{theorem}

        \begin{proof} If we prove $(i)$, then parts $(ii)$, $(iii)$
          and $(iv)$ will follow from $(i)$,
          Corollary~\ref{leaf-independence} and Corollary~\ref{basic}.
          
          We prove part $(i)$ by induction on the number of facets $q$
          of $\D$.  If $q>1$, then $q \geq 3$ (if $\D$ is a tree of
          two facets, both facets must be leaves by Lemma~\ref{two
          leaves}, and since $\D$ is connected, we can get minimal
          vertex covers of cardinalities one and two, which means that
          $\D$ is not unmixed).

          So the base case for induction is when $q=3$. In this case,
          let $F_1$ and $F_2$ be the two disjoint leaves of $\D$, and
          let $G_1$ be the third facet. Since $\D$ is connected and
          unmixed, $G_1$ cannot be a leaf (because the leaves are
          pairwise disjoint). So $G_1$ is a joint for both $F_1$ and
          $F_2$ and this settles the case $q=3$.

        For the general case, suppose that $G$ is a joint of $\Delta$. By
        Corollary~\ref{joint-not-leaf}, $G$ is not a leaf. By
        Corollary~\ref{unmixed-corollary} and
        Proposition~\ref{joint-removal-unmixed}, if we remove $G$, the
        forest $\D'=\D \rmv{G}$ is unmixed and
        $\al(\D')=r$. By the induction hypothesis,
        \begin{eqnarray}\label{DD'}\D'= \tuple{F_1, \ldots, F_r} \cup
        \tuple{G_1, \ldots, G_s}\end{eqnarray} where conditions $(i)$
        to $(iv)$ are satisfied. It is easy to see from condition $(iv)$ 
        that if $F$ is a
        leaf of $\D$, then it will still be a leaf of $\D'$ (because
        it has a free vertex). 

        Our goal is to show that the converse is true, that is, to
        show that $F_1, \ldots, F_r$ are all the leaves of $\D$.

        We have the following presentation for $\D$:
        \begin{eqnarray}\label{DD}\D= \tuple{F_1, \ldots, F_r} \cup
        \tuple{G_1, \ldots, G_s} \cup \tuple{G}.\end{eqnarray} There
        are two cases to consider.

        \begin{enumerate}

        \item[\emph{Case 1.}] $G$ is the only joint of $\D$. 

         Suppose, without loss of generality, that for some $e$, $F_1,
        \ldots, F_{e-1}$ are leaves of $\D$ and $F_e, \ldots, F_r$ are
        not leaves of $\D$. Remove $F_1, \ldots, F_{e-1}$ from $\D$ to
        obtain the forest $$\D''=\tuple{F_e, \ldots, F_r} \cup
        \tuple{G_1, \ldots, G_s} \cup \tuple{G}.$$ By Lemma~\ref{two
        leaves}, $\D''$ has at least two leaves. Neither one of $G_1,
        \ldots, G_s$ could be a leaf, because neither one of them has
        a free vertex. To see this, note that by the induction
        hypothesis on $\D'$ and part ($iv$) of the theorem, $G_1,
        \ldots, G_s$ do not have free vertices in $\D'$, and hence
        they cannot have free vertices in $\D$. As facets of $\D''$,
        they still do not have free vertices, because as $G$ is the
        only joint of $\D$, $$G_i \cap F_j \subseteq G \cap F_j
        \subseteq G{\rm \ \ for\ \ } 1 \leq i \leq s {\rm \ \ and\ \ }
        1 \leq j \leq e-1.$$ Since $G$ is a facet of $\D''$ the
        removal of $F_1, \ldots, F_{e-1}$ does not free any vertices
        of $G_1, \ldots, G_s$.

        This implies that at least one of $F_e, \ldots, F_r$ is a leaf
        of $\D''$. Suppose that $F_e$ is a leaf. Then there exists a
        facet $G' \in \D''$ such that
        $$H \cap F_e \subseteq G' \cap F_e \ \ {\rm for \ all \ } H
        \in \D''\rmv{F_e}.$$ Since $F_i \cap F_e = \emptyset$ for
        $i=1,\ldots e-1$, it follows that
        $$H \cap F_e \subseteq G' \cap F_e \ \ {\rm for \ all \ } H
        \in \D \rmv{F_e}$$ and so $F_e$ is a leaf of $\D$, which is a
        contradiction.

            \idiot{For $F_e$ it was picked to belong to a connected
          component of $\D''$ with at least two facets. We don't
          really need to worry here about disconnection; $F_e, \ldots,
          F_r$ don't really even get disconnected totally, as they
          never intersected $F_1, \ldots, F_{e-1}$ in the first
          place.}

        \item[\emph{Case 2.}] $\D$ has another joint $G'$
        distinct from $G$. 

        Consider the presentation of $\D$ as
        in~(\ref{DD}). As $\{F_1,\ldots, F_r\}$ is a maximal
        independent set of facets in $\D$, it cannot contain $G'$
        (Corollary~\ref{joint-not-leaf}). Therefore $$G' \in \{G_1,
        \ldots, G_s\}.$$
        
        We show that, say, $F_1$ is a leaf of $\D$. 

        Consider the two unmixed forests $$\D' =\D \rmv{G} {\rm \ \ and
        \ \ } \D'' =\D \rmv{G'}.$$

        We already know from before that $F_1$ is a leaf of $\D'$.
        From the fact that $\{F_1,\ldots, F_r\}$ is a maximal
        independent set of facets in $\D''$ and
        Corollary~\ref{leaf-independence} and the induction
        hypothesis, it follows that $F_1$ is also a leaf of $\D''$.
      
        So, by the definition of a leaf, there is a facet, say $G_1$,
        in $\D'$, such that
        \begin{eqnarray}\label{label-1} H \cap F_1 \subseteq G_1 
          \cap F_1 {\rm \ \ for \ all\ \ } H
        \neq G, F_1.\end{eqnarray}
      and a facet $G_2 \in \D''$ such that
      \begin{eqnarray}\label{label-2} H \cap F_1 \subseteq G_2 \cap
        F_1 {\rm \ \ for \ all\ \ } H \neq G', F_1.\end{eqnarray}

       The possible scenarios are the following.
   
      \begin{enumerate}   
      
       \item $G_1 \neq G'$ or $G_2 \neq G$.
      
       Suppose $G_1 \neq G'$. In this case $G_1 \in \D''$, and so
      because of~(\ref{label-2}) $$G_1 \cap F_1 \subseteq G_2 \cap
      F_1$$ which with~(\ref{label-1}) and~(\ref{label-2}) implies
      that
      $$H \cap F_1 \subseteq G_2 \cap F_1 {\rm \ \ for \ all\ \ } H
      \neq F_1$$
      hence $F_1$ is a leaf of $\D$.  The case $G_2 \neq G$ is identical.

      \item $G_1 = G'$ and $G_2 = G$.

        In this case, Statements (\ref{label-1}) and (\ref{label-2}),
        respectively, translate into
 
       \begin{eqnarray}\label{case21}\left \{ \begin{array}{ll} 
        H \cap F_1 \subseteq G' \cap F_1 & 
              {\rm \ \ for \ all\ \ } H \neq G,F_1 \\ 
        H \cap F_1  \subseteq G \cap F_1 &
              {\rm \ \ for \ all\ \ } H \neq G',F_1
	\end{array} 
       \right. \end{eqnarray}

        If $F_1$ is not a leaf of $\D$, it follows from (\ref{case21})
        that
        \begin{eqnarray}\label{case2}  \left \{ \begin{array}{l}
        G \cap F_1 \not \subseteq G' \cap F_1 \\
         G' \cap F_1 \not \subseteq G \cap F_1\\
         H \cap F_1 \subseteq
        (G\cap G') \cap F_1 {\rm \ \ for \ all\ \ } H \neq G, G', F_1 
        \end{array} \right. {}
        \end{eqnarray}
        
        \idiot{ This was not used in the proof. I just put it here: In
        the forest $\tuple{G, G', F_1}$ one can see that $G \cap G'
        \subseteq F_1.$ This is because $F_1$ is not a leaf, otherwise
        $F_1 \cap G \subseteq F_1 \cap G'$ or $F_1 \cap G' \subseteq
        F_1 \cap G$, which is a contradiction. So $G$ and $G'$ are
        leaves $\Longrightarrow$ $G \cap G' \subseteq G \cap F_1$, as
        $G \cap F_1 \subseteq G \cap G'$ is not possible.}

        By~(\ref{case2}) there exist \begin{eqnarray}\label{x-and-y}
           x \in (G \cap F_1) \setminus G'\
         \ {\rm and}\ \ y \in (G' \cap F_1) \setminus G.\end{eqnarray}

        {\bf Claim:} \emph{There is a minimal vertex cover for $\D
        \rmv{ G, G', F_1}$ that avoids all the vertices in $ G$, $G'$
        and $F_1$.}
        \vspace{.2in}

        \begin{itemize}

        \item[ ]\noindent \emph{Proof of Claim:} We first show that
        there is no facet of $\D\rmv{ G, G', F_1}$ that has all its
        vertices in $G\cup G'$. Suppose that $H$ is such a facet:
        \begin{eqnarray}\label{H} H = (H\cap G) \cup (H\cap
        G')\end{eqnarray} and consider the tree $$\D_1=
        \tuple{G,G',F_1,H}.$$ By Lemma~\ref{two leaves}, $\D_1$ must
        have two leaves. Note that $H$ cannot be a leaf, since because
        of (\ref{H}) it has no free vertices. If $F_1$ is a leaf, then
        it cannot have $G$ or $G'$ as its joint, since that violates
        the first two conditions in~(\ref{case2}), and so $H$ must be
        its joint. But then it follows that $$G \cap F_1 \subseteq H
        \cap F_1.$$ This implies that $x \in H$ (where $x$ is defined
        in (\ref{x-and-y})), which along with the third part
        of~(\ref{case2}), results in $x \in G'$, which is a
        contradiction.
        
        So $G$ and $G'$ are the two leaves of $\D_1$. Consider $G$
        first. If $G'$ is a joint for $G$, it follows that $$F_1 \cap
        G \subseteq G' \cap G \subseteq G'$$ which
        contradicts~(\ref{case2}).

        If $H$ is a joint of $G$, then
        $$F_1 \cap G \subseteq H \cap G $$ which implies that $x \in
        H$, but this again means $x \in G'$ (because of
        (\ref{case2})), which is a contradiction. So $F_1$ is the only
        possible joint for $G$.
        
        With an identical argument for $G'$, it follows that $F_1$ is
        a joint for both $G$ and $G'$ in $\D_1$, and therefore $$H
        \cap G \subseteq F_1 \cap G {\rm \ \ and \ \ } H \cap G'
        \subseteq F_1 \cap G'$$ which along with~(\ref{H}) implies
        that $$ H \subseteq F_1$$ which is impossible since $H$ and
        $F_1$ are both facets of $\D$.

        So we have shown that every facet of $\D$ other than $G$, $G'$
        and $F_1$, has at least one vertex outside $G$ and $G'$ (and
        therefore by the third condition in~(\ref{case2}), outside
        $F_1$).

        For each facet $H$ of $\D \rmv{ G, G', F_1}$, pick a vertex $z
        \in H$ that avoids all three facets $G$, $G'$ and $F_1$. The
        set of these vertices is a vertex cover for $\D \rmv{ G, G',
        F_1}$, and so it has a subset that is a minimal vertex cover.
        This proves the claim.
        \end{itemize}

        Now let $A$ be a minimal vertex cover for $\D \rmv{ G, G',
        F_1}$ that avoids all the vertices in $ G$, $G'$ and $F_1$.
        Since $\D \rmv{ G, G', F_1}$ has $r-1$ independent
        facets, $|A| \geq r-1$. Now $A \cup \{x,y\}$ is a minimal
        vertex cover for $\D$ with more than $r$ vertices, which
        contradicts the fact that $\D$ is unmixed with vertex covering
        number equal to $r$ (Note that $x$ and $y$ do not belong to
        any facet of $\D \rmv{ G, G', F_1}$, as this would
        contradict the third condition in~(\ref{case2})).
       \end{enumerate}
        \end{enumerate}

      So both cases 1 and 2 lead to contradictions, therefore all of
      $F_1, \ldots, F_r$ must be leaves of $\D$, which proves the
      theorem.
 \end{proof}

\begin{example}\label{structure-example} The simplicial complex $\D$ shown 
below is an unmixed tree, satisfying properties (i) to (iv) of
Theorem~\ref{structuretheorem}.


\[ \includegraphics{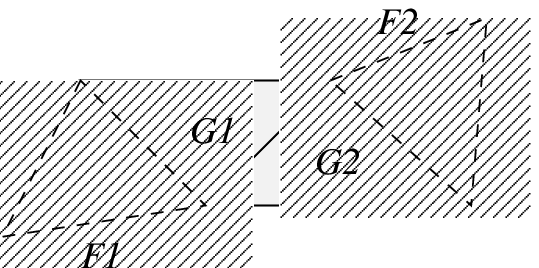} \]

\end{example}

It is important to notice that Theorem~\ref{structuretheorem} does not
suggest that every facet in an unmixed tree is either a leaf or a
joint (See Example~\ref{not-leaf-or-joint-example} below).  On the
other hand two different leaves in an unmixed tree may share a joint,
as is the case with the unmixed graph $\tuple{xy,yz,zu}$. For these
reasons the two numbers  $r$ and $s$ in the statement of
Theorem~\ref{structuretheorem} that count the number of leaves and
non-leaves, respectively, do not seem to have any particular relationship
to one another.

\begin{example}\label{not-leaf-or-joint-example} The following 
simplicial complex, which is the facet complex of the ideal
$$(xu,uvew,zvew,efw,efg,fgy)$$
is an unmixed tree with a facet
$\{e,f,w\}$ that is neither a leaf nor a joint. In fact, the two
leaves $\{ x,u\}$ and $\{z,v,e,w\}$ share a joint $\{u,v,e,w\}$.

  \[ \includegraphics{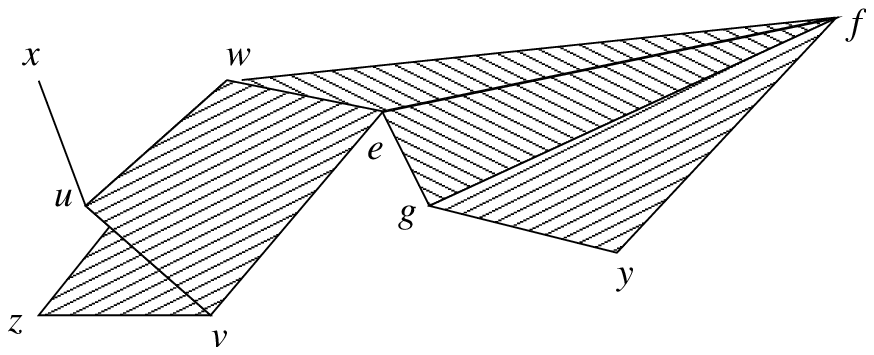} \] 

Above, for simplicity, an $n$-dimensional facet (simplex) is drawn as
a shaded polygon with $n+1$ vertices. The picture in 3D is as follows:

  \[ \includegraphics{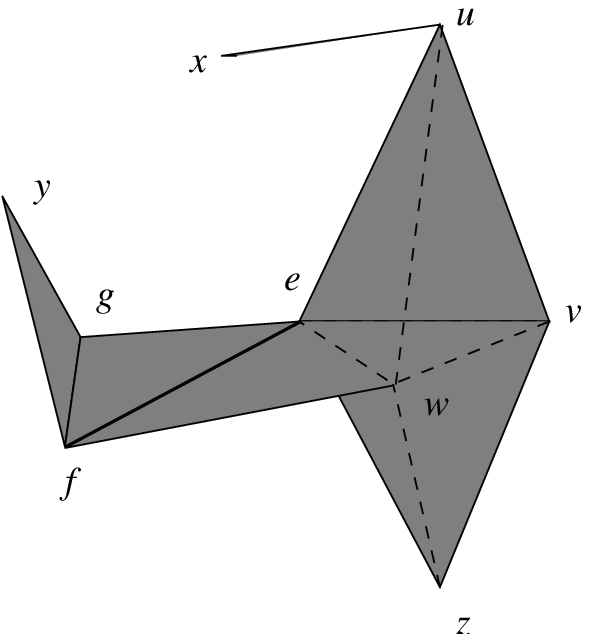} \]

\end{example}


\section{Grafting simplicial complexes}\label{grafting-section}

All that we proved in the previous section about unmixed trees can be
put into one definition-- namely that of a grafted tree. In fact, the
method of grafting works as an effective way to build an unmixed
simplicial complex from any given simplicial complex by adding new
leaves (Theorem~\ref{grafted-is-unmixed}). It turns out that a grafted
simplicial complex is Cohen-Macaulay (Theorem~\ref{CM}).

\begin{definition}[grafting]\label{grafting} A simplicial complex $\D$
  is a \emph{grafting} of the simplicial complex $\D'=\tuple{G_1,
    \ldots, G_s}$ with the simplices $F_1, \ldots, F_r$ (or we say
    that $\D$ is \emph{grafted}) if $$\D= \tuple{F_1, \ldots, F_r}
    \cup \tuple{G_1, \ldots, G_s}$$ with the following properties:

\begin{enumerate}
\item[(i)] $V(\D') \subseteq V(F_1) \cup  \ldots \cup V(F_r)$;
\item[(ii)] $F_1, \ldots, F_r$ are all the leaves of $\D$;
\item[(iii)] $\{G_1, \ldots, G_s\} \cap \{F_1, \ldots, F_r\} = \emptyset$;
\item[(iv)] For $i \neq j$, $F_i \cap F_j = \emptyset$;
\item[(v)] If $G_i$ is a joint of $\Delta$, then $\D \rmv{G_i}$ is
  also grafted.

\end{enumerate} 
\end{definition}

Note that a simplicial complex that consists of only one facet or
several pairwise disjoint facets is indeed grafted, as it could be
considered as a grafting of the empty simplicial complex. It is easy
to check that conditions (i) to (v) above are satisfied in this case.

It is also clear that the union of two or more grafted simplicial
complexes is itself grafted.

\begin{remark}\label{embedded-intersections} Condition~(v)
  above implies that if $F$ is a leaf of a grafted $\D$, then all the
  facets $H$ that intersect $F$ have embedded intersections; in other words
  if $H \cap F$ and $H' \cap F$ are both nonempty, then 
  $$H \cap F \subseteq H' \cap F \ \ {\rm or}\ \ H' \cap F \subseteq H
  \cap F.$$
  
  This implies that there is a chain of intersections $$H_1 \cap F
  \supseteq \ldots \supseteq H_t \cap F$$
  where $H_1, \ldots, H_t$ are
  all the facets of $\D$ that intersect $F$.
\end{remark}

\begin{remark}\label{remove-all-grafted} Condition (v) in
  Definition~\ref{grafting} can be replaced by ``$\D \rmv{G_i}$ is
  grafted for all $i=1,\ldots, s$''. This is because even if $G_i$ is
  not a joint of $\D$, $\D \rmv{G_i}$ satisfies properties (i), (iii)
  and (iv), and it satisfies (ii) and (v) because of
  Remark~\ref{embedded-intersections}, and so $\D \rmv{G_i}$ is
  grafted. \end{remark}

\begin{remark}[A grafting of a tree is also  a tree] If $\D'$ in 
  Definition~\ref{grafting} is a tree, then $\D$ is also a tree. To
  see this, consider any subcollection $\D''$ of $\D$. If $\D''$
  contains $F_i$ for some $i$, then by
  remarks~\ref{embedded-intersections} and~\ref{remove-all-grafted}
  $F_i$ is a leaf of $\D''$. If $\D''$ contains neither of the $F_i$,
  then it is a subcollection of the tree $\D'$, which implies that
  $\D''$ has a leaf.
\end{remark}

The ``suspension'' of a graph, as defined in~\cite{Vi1}, is also a grafting
of that graph.

\begin{example}\label{grafting-example} The tree 
  $\tuple{F_1,F_2,G_1,G_2}$ that appeared in
  Example~\ref{structure-example} above is a grafting of the tree
  $\tuple{G_1,G_2}$ with the leaves $F_1$ and $F_2$.  In fact, there
  may be more than one way to graft a given simplicial complex.  For
  example, some possible ways of grafting $\tuple{G_1,G_2}$ are shown
  below:

\begin{picture}(300,220)
\put(0,180){$\D:$}
\put(40,155){\epsfig{file=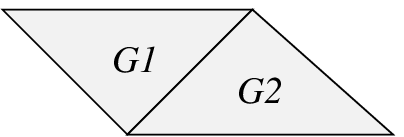, height=.5in}}
\put(150,180){$\stackrel{\mbox{\emph{graft}}}{\vector(1,0){50}}$}
\put(220,180){$\D':$}
\put(250,150){\epsfig{file=basic-grafted.eps, height=1in}}
\put(100,130){$\vector(0,-1){50}$}
\put(110,105){\emph{graft}}
\put(0,30){$\D'':$} 
\put(40,0){\epsfig{file=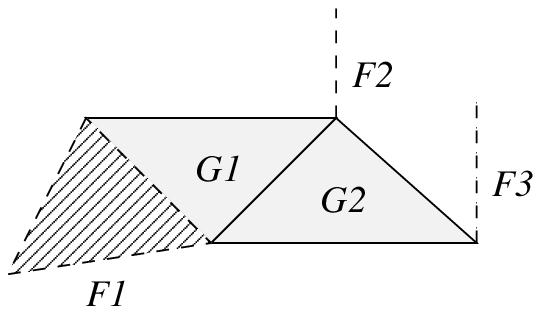, height=1.1in}}
\put(180,130){$\vector(1,-1){50}$}
\put(220,105){\emph{graft}}
\put(220,30){$\D''':$}
\put(250,10){\epsfig{file=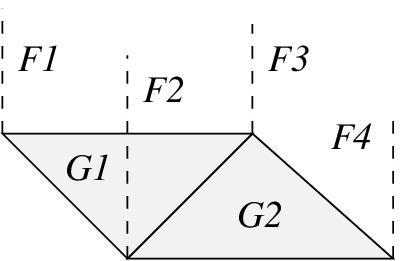, height=1in}}
\end{picture}

\end{example}


\begin{theorem}[A grafted simplicial complex is 
  unmixed]\label{grafted-is-unmixed} Let $$\D= \tuple{F_1, \ldots,
    F_r} \cup \tuple{G_1, \ldots, G_s}$$ be a grafting of the
    simplicial complex $\tuple{G_1, \ldots, G_s}$ with the simplices
    $F_1, \ldots, F_r$. Then $\D$ is unmixed, and $\alpha(\D)=r$.
\end{theorem}

       \begin{proof}  If $\tuple{G_1, \ldots, G_s}$ is the 
         empty simplicial complex, the claim is immediate, so we
         assume that it is nonempty. 

         We argue by induction on the number of facets $q$ of
         $\D$. The first case to consider is $q=3$. In this case, $\D$
         must have at least two leaves, as if there were only one leaf
         $F_1$, i.e. if $\D = \tuple{F_1} \cup \tuple{G_1, G_2}$, then
         by Condition (i) of Definition~\ref{grafting} we would have
         $G_1 \subseteq F_1$ and $G_2 \subseteq F_1$, which is
         impossible. So $\D = \tuple{F_1, F_2} \cup \tuple{G_1}$,
         where $G_1 \subseteq F_1 \cup F_2$ and $F_1 \cap F_2
         =\emptyset$. It is now easy to see that $\D$ is unmixed with
         $\al(\D)=2$.

         Suppose $\D$ has $q >3$ facets, and let $G_1$ be a joint of
         the leaf $F_1$. By Part (v) of Definition~\ref{grafting} $\D'
         = \D \rmv{G_1}$ is also grafted, and therefore by the
         induction hypothesis unmixed with $\al(\D')=r$.
         
         Let $A$ be a minimal vertex cover of $\D$. We already know
         that $|A| \geq r$ as $F_1, \ldots, F_r$ are $r$ independent
         facets of $\D$. Now suppose that $|A| > r$ . Since $A$ is
         also a vertex cover for $\D'$, it has a subset $A'$ that is a
         minimal vertex cover of $\D'$ with $|A'|=r$. Since $A'$ is a
         proper subset of $A$, it is not a vertex cover for $\D$, and
         therefore $A'$ cannot contain a vertex of $G_1$. So $A'$
         contains a free vertex $x$ of $F_1$ (all non-free vertices of
         $F_1$ are shared with $G_1$).  Now $A$ must contain a vertex
         $y$ of $G_1$; say $y \in G_1 \cap F_2$ ($y \notin F_1$, since
         in that case $x$ would be redundant). So $$A=A' \cup \{y\}.$$
         On the other hand $A'$ must also contain a vertex of $F_2$,
         say $z$. So $F_2$ contributes two vertices $y$ and $z$ to
         $A$; note that neither one of $y$ or $z$ could be a free
         vertex, as in that case the free one would be redundant.

         Now suppose that $G_2$ is a joint of $F_2$. Remove $G_2$ from
         $\D$ to get $$\D'' = \D \rmv{G_2}.$$
         So $A$ has a subset
         $A''$, $|A''|=r$, that is a minimal vertex cover for $\D''$.
         But as $A$ already has exactly one vertex in each of $F_1,
         F_3, \ldots, F_r$, the only way to get $A''$ from $A$ is to
         remove one of $y$ or $z$, this means that:
$$A'' = A\setminus \{y\} \ \ \ \ \ {\rm{or}}\ \ \ \ \ A'' = A\setminus
\{z\}.$$ In either case $A''$  contains a vertex of $G_2$, which
implies that $A''$ is a minimal vertex cover for $\D$; a contradiction.
\end{proof}

Example~\ref{grafting-example} demonstrates
Theorem~\ref{grafted-is-unmixed}: $\D=\tuple{G_1,G_2}$ is a
non-unmixed tree, which gets grafted with some leaves to make the
unmixed trees $\D'$, $\D''$ and $\D'''$.

One could graft any simplicial complex, even a badly non-unmixed non-tree.

\begin{example} Let $\D'$ be the non-unmixed non-tree in 
  Example~\ref{non-leaf-example}. We could graft $\D'$ with three new
  leaves $$\{x,y,v\}, \{u,w\}, \{z,e\}$$
    
  The resulting picture below is unmixed, and moreover, as we prove
  later, Cohen-Macaulay.
\[ \includegraphics[height=1in]{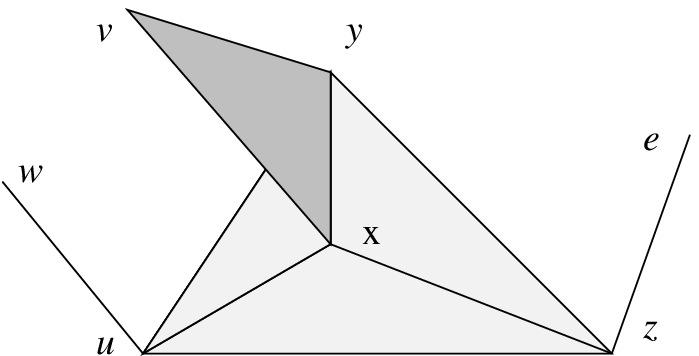} \]
\end{example}

In the case of a tree theorems~\ref{structuretheorem}
and~\ref{grafted-is-unmixed} put together with Corollary~\ref{basic}
produce a much stronger statement:

\begin{corollary}[A tree is unmixed if and only if 
grafted]\label{tree-unmixed-iff-grafted} Suppose the simplicial
  complex $\D$ is a tree. Then $\D$ is unmixed if and only if $\D$ is
  grafted. \end{corollary}

Grafted simplicial complexes behave well under localization; in other
words, the localization of a grafted simplicial complex is also
grafted. In the case of trees this follows directly from
Corollary~\ref{tree-unmixed-iff-grafted},
Lemma~\ref{unmixed-localizes} and Lemma~\ref{localization}.  But the
statement holds more generally.

\begin{proposition}[Localization of a grafted simplicial complex is 
grafted]\label{grafted-localization} Let $I=\F(\D)$ where $\D$ is a
  grafted simplicial complex with vertices labeled $\xs$. Suppose that
  $k$ is a field and $p$ is a prime ideal of the polynomial ring
  $k[\xs]$. Then $\df(I_p)$ is a grafted simplicial complex.
\end{proposition}

           \begin{proof} With notation as in Definition~\ref{grafting}, 
             let $$\D= \tuple{F_1, \ldots, F_r} \cup \tuple{G_1,
               \ldots, G_s}.$$
             
             If $\D$ has only one facet, the statement of the theorem
             follows immediately, so assume that $\D$ has two or more
             facets.
             
             As in the proof of Lemma~\ref{localization}, it is enough
             to assume that $p$ is generated by a subset of $\{\xs
             \}$, so
             $$p=(x_{i_1},\ldots, x_{i_h}).$$

             Following Discussion~\ref{local-discussion}, let
             $$\D_p=\df(I_p)=\tuple{F'_1, \ldots, F'_t} \cup \tuple{G'_1,
               \ldots, G'_u}$$ where for $i=1,\ldots,r$ and $j=1,\ldots,s$
             $$F'_i=F_i \cap \{ x_{i_1},\ldots, x_{i_h}\} \ \ {\rm
               and}\ \ G'_j=G_j \cap \{ x_{i_1},\ldots, x_{i_h}\}$$
             and $F'_{t+1}, \ldots, F'_r, G'_{u+1}, \ldots, G'_s$ each
             contain at least one of
             \begin{eqnarray}\label{the-extra-facets}F'_1, \ldots,
               F'_t, G'_1, \ldots, G'_u. \end{eqnarray}
             
             We now rename the facets of $\D_p$ as follows. For
             $i=1,\ldots,t$, let $$H_i=F'_i.$$
             
             For each $i=t+1,\ldots, r$, $F'_i$ contains one of the
             facets appearing in~(\ref{the-extra-facets}). But as by
             definition $F_i \cap F_j = \emptyset$ for all $j \neq i$,
             there must be some $j \leq u$ for which $G'_j \subseteq
             F'_i$. For this particular $j$, set
             $$H_i=G'_j.$$
             
             This choice of $j$ is well-defined: if there were some $f
             \leq u$ distinct from $j$ such that $G'_f \subseteq
             F'_i$, then it would follow from
             Remark~\ref{embedded-intersections} that either $G'_j
             \subseteq G'_f$ or $G'_f \subseteq G'_j$, which
             contradicts the fact that both $G'_j$ and $G'_f$ are
             facets of $\D_p$.

             We now represent $\D_p$ as $$
             \D_p=\tuple{H_1,\ldots,H_r}
             \cup \tuple{E_1,\ldots,E_v}$$
             where $E_1,\ldots,E_v$
             represent all the other facets of $\D_p$ that were not
             labeled by some $H_i$.
             
             Our goal is to show that $\D_p$ is a grafting of the
             simplicial complex $\tuple{E_1,\ldots,E_v}$ with the
             simplices $H_1,\ldots,H_r$.

             It is clear by our construction that the facets
             $H_1,\ldots,H_r$ are pairwise disjoint. To see this,
             notice that for each pair of distinct numbers $i_1, i_2
             \leq r$, there is a pair of distinct numbers $j_1, j_2
             \leq r$ such that $$H_{i_1} \subseteq F'_{j_1}\subseteq
             F_{j_1}\ \ {\rm and}\ \ H_{i_2} \subseteq
             F'_{j_2}\subseteq F_{j_2}$$ and as $F_{j_1} \cap F_{j_2}
             =\emptyset$,
             $$H_{i_1} \cap H_{i_2} = \emptyset.$$
             
             So Condition~(iv) of Definition~\ref{grafting} is
             satisfied.
             
             On the other hand, by Theorem~\ref{grafted-is-unmixed}
             $\D$ is unmixed, so by Lemma~\ref{unmixed-localizes}
             $\D_p$ is unmixed with $\al(\D_p)=\al(\D)=r$. We now
             apply Lemma~\ref{basic-lemma} to $\D_p$ to deduce that
             $$V(\D_p)=V(H_1)\cup \ldots \cup V(H_r),$$
             which implies
             Condition~(i) in Definition~\ref{grafting}. This also
             implies that $E_1,\ldots,E_v$ cannot have free vertices,
             and hence cannot be leaves of $\D_p$.
             
             Condition~(iii) is satisfied by the construction of
             $\D_p$.
             
             We need to show that $H_1,\ldots,H_r$ are all leaves of
             $\D_p$. If $\D_p=\tuple{H_1,\ldots,H_r}$ then $\D_p$ is
             grafted by definition. So suppose that $\D_p$ has a
             connected component $\D'$ with two or more facets. As
             $\D'$ is connected, it must contain some of the $E_i$,
             and as $V(\D_p)=V(H_1)\cup \ldots \cup V(H_r)$, $\D'$
             must also contain some of the $H_j$. So we can without
             loss of generality assume that
             $$\D'= \tuple{H_1,\ldots,H_e} \cup
             \tuple{E_1,\ldots,E_f}$$
             for some $1 \leq e \leq r$ and $1
             \leq f \leq v$.
             
             We now show that, for example, $H_1$ is a leaf for $\D'$.
             There are two cases to consider:

             \begin{enumerate}
             \item[\emph{Case 1.}] $H_1=F'_i$ for some $i$ such that
               $1\leq i \leq t$. 
               
               Since $\D'$ is connected, it has some facets that
               intersect $H_i$; suppose that $E_{j_1}, \ldots,
               E_{j_l}$ are all the facets of $\D'\rmv{H_1}$ such that
               $$H_1 \cap E_{j_z} \neq \emptyset$$ for $z=1,\ldots,l$.
               
               For each $z=1,\ldots,l$ suppose that $$E_{j_z} =
               G'_{m_z}.$$
               
               The above paragraph translates into
               $$F'_i \cap G'_{m_z} \neq \emptyset$$
               and hence
               $$F_i\cap G_{m_z} \neq \emptyset$$
               for $z=1,\ldots,l$.
               
               From Remark~\ref{embedded-intersections} it follows
               that there is some total order of inclusion on the
               nonempty sets $F_i \cap G_{m_z}$; we assume that $$F_i
               \cap G_{m_1} \supseteq F_i \cap G_{m_2} \supseteq
               \ldots \supseteq F_i \cap G_{m_l}$$
               which after
               intersecting each set with $\{ x_{i_1},\ldots,
               x_{i_h}\}$ turns into
               $$F'_i \cap G'_{m_1} \supseteq F'_i \cap G'_{m_2} \supseteq 
               \ldots \supseteq F'_i \cap G'_{m_l}$$
               which is equivalent to 
               $$H_1 \cap E_{j_1} \supseteq H_1 \cap E_{j_2} \supseteq
               \ldots \supseteq H_1 \cap E_{j_l}$$
               
               It follows that $H_1$ is a leaf of $\D'$, and in
               addition, Condition~(v) of Definition~\ref{grafting} is
               satisfied.

             \item[\emph{Case 2.}] $H_1=G'_j$ for some $j$ such that
               $1\leq j \leq u$.
               
               In this case for some $i$, $t < i \leq r$, $$
               H_1=G'_j\subseteq F'_i.$$
               
               Exactly as above, let $E_{j_1}, \ldots, E_{j_l}$ be all
               the facets of $\D'\rmv{H_1}$ such that $H_1 \cap
               E_{j_z} \neq \emptyset$, and let $E_{j_z} = G'_{m_z}$
               for $z=1,\ldots,l$.  
               
               As all the sets $F_i \cap G_{m_z}$ are nonempty, we
               follow the exact argument as above to obtain the chain
               $$F'_i \cap G'_{m_1} \supseteq F'_i \cap G'_{m_2}
               \supseteq \ldots \supseteq F'_i \cap G'_{m_l}$$
               
               As $G'_j\subseteq F'_i$, we can intersect all these
               sets with $G'_j$ to obtain 
               $$G'_j \cap G'_{m_1} \supseteq G'_j \cap G'_{m_2}
               \supseteq \ldots \supseteq G'_j  \cap G'_{m_l}$$
               which is equivalent to
               $$H_1 \cap E_{j_1} \supseteq H_1 \cap E_{j_2}
               \supseteq \ldots \supseteq H_1  \cap E_{j_l}$$
               
               It follows that $H_1$ is a leaf of $\D'$, and also
               Condition~(v) of Definition~\ref{grafting} is
               satisfied.
               \end{enumerate}
               
           \end{proof}

        
\section{Grafted simplicial complexes are Cohen-Macaulay}\label{CM-section}

We are now ready to show that the facet ideal of a grafted simplicial
complex has a Cohen-Macaulay quotient. Besides revealing a wealth of
square-free monomial ideals with Cohen-Macaulay quotients, this result
implies that all unmixed trees are Cohen-Macaulay.


Let $\D$ be a grafted simplicial complex over a vertex set
$V=\{\xs\}$. By Definition~\ref{grafting}, $\D$ will have the
form $$\D= \tuple{F_1, \ldots, F_r} \cup \tuple{G_1, \ldots, G_s}$$
where $\al(\D)=r$ and $F_1, \ldots, F_r$ are the leaves of $\D$.

Let $$\R(\D)= k[\xs]/\F(\D),$$ where $k$ is a field and let
$$\m=(\xs)$$ be the irrelevant maximal ideal.

From Discussion~\ref{dimension} we know that $$\dimn \R(\D)=n-r.$$

In order to show that $\R(\D)$ is Cohen-Macaulay, it is enough to show
that there is a homogeneous regular sequence in $\m$ of length $n-r$.

\idiot{See Corollary 1.1.3, and 1.5.8 and 1.5.9 of~\cite{BH}.}

It is interesting to observe how the previous sentence follows also
from Proposition~\ref{grafted-localization}: if $m$ is any other
maximal ideal of $\R(\D)$, from the proof of
Lemma~\ref{localization} and Proposition~\ref{grafted-localization} we
see that if $p=(x_1,\ldots, x_e)$ is the ideal generated by all of
$x_i$ that belong to $m$, then $I_{m}=I_p$ is the facet ideal of a
grafted simplicial complex over the vertex set $\{x_1,\ldots,
x_e\}$. So one can write $m=p+q$ where $q$ is a prime ideal of
$k[x_{e+1},\ldots, x_n]$. It follows that
$$\R(\D)_{m}=k[x_1,\ldots,x_e]_p/I_p \otimes_k
k[x_{e+1},\ldots, x_n]_q.$$

As $k[x_{e+1},\ldots, x_n]_q$ is clearly Cohen-Macaulay, by
Theorem~5.5.5 of~\cite{V}, it is enough to show that
$k[x_1,\ldots,x_e]_p/I_p$ is Cohen-Macaulay in order to conclude that
$\R(\D)_{m}$ is Cohen-Macaulay.  But this is again the case of
localizing at the irrelevant ideal.

Now suppose that for each $i \leq r$, $$F_i=y_i x^i_1\ldots x^i_{u_i}$$
where $y_i$ is a free vertex of the leaf $F_i$, and $y_i, x^i_1,
\ldots, x^i_{u_i} \in V$. We wish to show that 
\begin{eqnarray}\label{the-sequence} y_1 - x^1_1, \ldots,
y_1 - x^1_{u_1}, \ldots, y_r - x^r_1, \ldots, y_r - x^r_{u_r}
 \end{eqnarray} is a regular sequence in $\R(\D)$. This follows from the
 process of ``polarization'' that we describe below.

\begin{proposition}[{\cite{Fr}}]\label{polarization} Let $R$ be the ring 
  $k[\xs]/(\Ms)$, where $\Ms$ are monomials in the variables $\xs$,
  and $k$ is a field.  Then there is an $N \geq n$, and a set of
  square-free monomials $N_1, \ldots, N_q$ in the polynomial ring
  $k[x_1, \ldots, x_N]$, such that
$$R=R'/(f_1,\ldots,f_{N-n})$$ where $R'=k[x_1, \ldots, x_N]/(N_1,
\ldots, N_q)$ and $f_1,\ldots,f_{N-n}$ is a regular sequence of forms
of degree one in $R'$.
\end{proposition}

For the purpose of our argument, it is instructive to see an
outline of the proof of this proposition.

        \begin{proof}[Sketch of proof] Suppose, without loss of
        generality, that $x_1 | M_i$ for $1 \leq i \leq s$, and $x_1
        \ndiv M_j$ for $s < j \leq q$.

        For $i=1, \ldots, s$ we set $$M_i' =\frac{M_i}{x_1}$$ so that
        we can write $$I=(\Ms)= (x_1M_1', \ldots,x_1M_s', M_{s+1},
        \ldots, M_q).$$

        Define $$I_1=(x_{n+1}M_1', \ldots,x_{n+1}M_s', M_{s+1},
        \ldots, M_q) \subseteq k[\xs,x_{n+1}].$$
        
        Then $R=R_1/(x_{n+1}-x_1)$ where $$R_1=k[\xs,x_{n+1}]/I_1.$$

        It is then shown that $x_{n+1}-x_1$ is a non-zerodivisor in
        $R_1$. If $I_1$ is square-free, we are done. Otherwise one
        applies the same procedure to $I_1$ continually until the
        ideal becomes square-free.  \end{proof}
        

What we would like to show is that Sequence~(\ref{the-sequence})
polarizes the ring $$S=k[y_1,\ldots,y_r]/(y_1^{u_1+1}, \ldots,
y_r^{u_r+1},E_1, \ldots, E_s)$$ into the ring $\R(\D)$, where $E_1,
\ldots, E_s$ are monomials corresponding to the facets $G_1, \ldots,
G_s$, where each vertex belonging to $F_i$ has been replaced by the
free vertex $y_i$. In other words if $$J=(y_1 - x^1_1, \ldots, y_1 -
x^1_{u_1}, \ldots, y_r - x^r_1, \ldots, y_r - x^r_{u_r}),$$ we wish to
show that $$S=\R(\D)/J.$$ It will then follow from the proof of
Proposition~\ref{polarization} (as detailed in~\cite{Fr} as well as
in~\cite{Vi2}) that Sequence~(\ref{the-sequence}) is a regular
sequence in $\R(\D)$.

Intuitively our claim is straightforward to see. The only problem that
may arise is if after applying Sequence~(\ref{the-sequence}) to $S$,
we end up with a permutation of the vertices of $\D$. To prevent this
from happening, we use the subtle structure of a grafted simplicial
complex (Remark~\ref{embedded-intersections}) that the facets
intersecting a leaf do so in an embedded (and therefore ordered)
manner. In other words, suppose for the leaf $F_i$, the facets $H^i_1,
\ldots, H^i_{e_i}$ are all the facets of $\D \rmv{F_i}$ that intersect
$F_i$, with the ordering
\begin{eqnarray}\label{ordered-sequence} H^i_1 \cap F_i \subseteq
\ldots \subseteq H^i_{e_i} \cap F_i.\end{eqnarray} So in
Sequence~(\ref{the-sequence}), we order
\begin{eqnarray}\label{sub-polar}y_i - x^i_1, \ldots, y_i -
x^i_{u_i}\end{eqnarray} such that if for any $e$ and $f$, $x^i_e \in H^i_f$ 
then $x^i_ e \in H^i_{f+1}$.

We now use induction on the number of facets of $\D$. If we remove a
joint, say $G_1 \in \U_{\D}(F_1)$, we obtain a grafted simplicial
complex $$\D'=\D \rmv{G_1}$$
over the same set of vertices $\xs$, with
$\al(\D') = \al(\D)$ (Lemma~\ref{unmixed-induction}). Therefore if
$$\R(\D') = k[\xs]/\F(\D') $$ then $$\dimn \R(\D)= \dimn \R(\D').$$

Moreover, $\D'$ has $F_1, \ldots, F_r$ as leaves. So by the induction
hypothesis, Sequence~(\ref{the-sequence}) polarizes the ring
$$S'=k[y_1,\ldots,y_r]/(y_1^{u_1+1}, \ldots, y_r^{u_r+1},E_2, \ldots,
E_s)$$ into $\R(\D')$, or in other words, $$S'=\R(\D')/J.$$

The induction hypothesis has ensured that after applying
Sequence~(\ref{the-sequence}) to $S'$, all facets of $\D'$ are
restored to their original positions and labeling. Now it all reduces
to showing that during this polarization process, $E_1$ turns into
$G_1$.

This is clear, as for every $i$, $G_1 \cap F_i$ has its place in the
ordered sequence~(\ref{ordered-sequence}), and so if $|G_1 \cap
F_i|=h_i$, then the first $h_i$ applications of
Sequence~(\ref{sub-polar}) restore $G_1 \cap F_i$ before moving on to
facets that have larger intersections with $F_i$. As $G_1$ has
disjoint intersections with $F_1, \ldots, F_r$, once
Sequence~(\ref{sub-polar}) has been applied for all $i$, $G_1$ is
restored to its proper position.

We have shown that:

\begin{theorem}[Grafted simplicial complexes are Cohen-Macaulay]\label{CM} 
Let $\D$ be a grafted simplicial complex over a set of vertices
labeled $\xs$, and let $k$ be a field. Then $\R(\D)= k[\xs]/\F(\D)$ is
Cohen-Macaulay.\end{theorem}

Theorem~\ref{CM} along with Proposition~\ref{CM-is-unmixed} and
Corollary~\ref{tree-unmixed-iff-grafted} imply that for a tree being
unmixed and being Cohen-Macaulay are equivalent conditions.

\begin{corollary}[A tree is Cohen-Macaulay if and only if 
unmixed]\label{CM-criterion} Let $\D$ be a tree over a set of
 vertices $\xs$, and let $k$ be a field. Then the quotient ring
 $k[\xs]/\F(\D)$ is Cohen-Macaulay if and only if $\D$ is
 unmixed. \end{corollary}


\end{document}